%% file: DCOPF.tex
\documentclass[journal,twoside,web]{ieeecolor}
\usepackage{generic}
\usepackage{cite}
\usepackage{amsmath,amssymb}
\usepackage{algorithmic}
\usepackage{graphicx}
\usepackage{textcomp}
\usepackage[intoc]{nomencl}
\makenomenclature
\usepackage{ifthen}
\renewcommand{\nomgroup}[1]{%
\ifthenelse{\equal{#1}{P}}{\item[\textbf{Parameters}]}{%
\ifthenelse{\equal{#1}{V}}{\item[\textbf{Variables}]}{%
\ifthenelse{\equal{#1}{S}}{\item[\textbf{sets}]}{%
\ifthenelse{\equal{#1}{O}}{\item[\textbf{Operators}]}{}}}}
}
\DeclareMathOperator{\diag}{diag}

\usepackage{accents}
\newcommand{\ubar}[1]{\underaccent{\bar}{#1}}

\usepackage{algorithm}

\newtheorem{theorem}{Theorem}
\newtheorem{prop}{Proposition}
\newtheorem{definition}{Definition}
\newtheorem{rem}{Remark}
\newtheorem{lem}{Lemma}

\usepackage{balance}

\usepackage{hyperref}

\allowdisplaybreaks

\DeclareSymbolFont{bbold}{U}{bbold}{m}{n}
\DeclareSymbolFontAlphabet{\mathbbold}{bbold}

\def\BibTeX{{\rm B\kern-.05em{\sc i\kern-.025em b}\kern-.08em
    T\kern-.1667em\lower.7ex\hbox{E}\kern-.125emX}}
\begin{document}
\title{Optimal Power Flow in DC Networks with Robust Feasibility and Stability Guarantees}
\author{Jianzhe~Liu,~\IEEEmembership{Member,~IEEE}, Bai~Cui,~\IEEEmembership{Member,~IEEE}, Daniel~K.~Molzahn,~\IEEEmembership{Senior Member,~IEEE}, Chen~Chen,~\IEEEmembership{Senior Member,~IEEE},
	Xiaonan~Lu,~\IEEEmembership{Member,~IEEE}, and Feng~Qiu,~\IEEEmembership{Senior Member,~IEEE}
\thanks{J. Liu and F. Qiu are with the Energy Systems Division, Argonne National Laboratory, Lemont, IL 60439, USA. Email: \{jianzhe.liu, fqiu\}@anl.gov. Argonne National Laboratory's work is based upon work supported by the U.S. Department of Energy’s Office of Energy Efficiency and Renewable Energy (EERE) under the Solar Energy Technologies Office Award Number 34230.}
\thanks{B. Cui is with the Power Systems Engineering Center, National Renewable Energy Laboratory, Golden, CO 80401, USA. Email: bcui@nrel.gov.}
\thanks{D. Molzahn is with the School of Electrical and Computer Engineering, Georgia Institute of Technology, Atlanta, GA 30332, USA. Email: molzahn@gatech.edu.}
\thanks{C. Chen is with the School of Electrical Engineering, Xi'an Jiaotong University, Shaanxi, China. Email: morningchen@xjtu.edu.cn.}
\thanks{X. Lu is with the College of Engineering, Temple University, Philadelphia, PA 19122, USA. Email: xiaonan.lu@temple.edu.}}

\maketitle

\begin{abstract}
With high penetrations of renewable generation and variable loads, there is significant uncertainty associated with power flows in DC networks such that stability and operational constraint satisfaction are of concern. Most existing DC network optimal power flow (DN-OPF) formulations assume exact knowledge of loading conditions and do not provide stability guarantees. In contrast, this paper studies a DN-OPF formulation which considers both stability and operational constraint satisfaction under uncertainty. The need to account for a range of uncertainty realizations in this paper's robust optimization formulation results in a challenging semi-infinite program (SIP). The proposed solution algorithm reformulates this SIP into a computationally tractable problem by constructing a tight convex inner approximation of the stability set using sufficient conditions for the existence of a feasible and stable power flow solution. Optimal generator setpoints are obtained by optimizing over the proposed convex stability set. The validity and value of the proposed algorithm are demonstrated through various DC networks adapted from IEEE test cases.
\end{abstract}

\input{nomenclature.tex}
\section{Introduction}
\label{sec:introduction}
Recent years have witnessed the growth of DC loads and generators, such as DC fast charging facilities, photovoltaic generation, and various electronic devices in sites like data centers. Interconnecting DC components in a DC network is an efficient operation method due to the reduction of DC-AC conversion stages~%
\cite{dorfler2018electrical}. DC networks have thus found promising applications in low- and medium-voltage power systems such as community nanogrids and microgrids, shipboard power systems, data centers, etc.~%
\cite{dragivcevic2016dc}. Common features of DC networks include: 1) the uncertainty in loading conditions is usually more significant due to the higher penetration level of uncertain components and relatively smaller scale compared to AC networks, and 2) many loads are controlled as constant power loads (CPLs) that have a destabilizing negative resistance effect that reduces the damping in a system~\cite{dragivcevic2016dc}.\footnote{DC CPLs can be considered as a special case of AC CPLs when the CPLs have a unity power factor. With a non-unity power factor, AC CPLs behave as negative resistances and reactances~\cite{hossain2018comprehensive}.}

Similar to other power systems, a DC network should work at a stable operating point that satisfies all operational constraints. The classic method to compute such an operating point is to formulate and solve an optimal power flow (OPF) problem. An OPF problem finds the optimal generation schedule corresponding to the system operating point that maximizes economic welfare while satisfying various physical and operational constraints.

Solving OPF problems to optimality with an AC power flow model is generally challenging due to the associated nonconvexity~\cite{lehmann2016ac}. Many research efforts have been devoted to improve OPF tractability using approximation and relaxation methods~\cite{molzahn2018survey}. Recent research has also studied OPF problems for DC networks (\mbox{DN-OPF})~%
\cite{farasat2015ga, gan2014optimal,li2018optimal,garces2018newton, inam2016stability,7464840}. Note that DN-OPF problems are fundamentally different from the so-called ``\mbox{DC-OPF}'' problems~\cite{stott2009dc}:
1) A DC-OPF problem is a simplified OPF problem for an AC system where the nonlinear AC power flow equations are linearized. Conversely, a \mbox{DN-OPF} problem considers the nonlinear power flow equations that accurately model the physics of a DC network. 2) While a DC-OPF problem is usually convex, a general DN-OPF problem is a nonconvex optimization problem~%
\cite{garces2018newton}. 

A variety of methods have been applied to solve \mbox{DN-OPF} problems. In~%
\cite{farasat2015ga}, a genetic algorithm is applied to solve the OPF problem for a DC distribution system. In~%
\cite{montoya2018linear}, linearization techniques are used to simplify the problem. Other methods~%
\cite{gan2014optimal, li2018optimal} employ second-order cone programming (SOCP) and quadratic convex programming to relax a DN-OPF problem into a convex formulation. The existing works demonstrate the capability to effectively solve various deterministic DN-OPF problems.

Despite recent advances, existing DN-OPF works in the literature have limitations in providing stability and feasibility guarantees when significant uncertainties are present. 
    First, previous results primarily focus on deterministic \mbox{DN-OPF} problems where the loading conditions are assumed to be fixed and known \textit{a priori}. Nevertheless, with high penetrations of intermittent generation and variable loads, uncertainty in the net loading conditions is a salient feature of DC networks~%
    \cite{liu2018robust}. Directly applying the OPF decisions computed using a specific scenario to an uncertain system can cause unpredictable deviations of the system operating point from the designated value~%
    \cite{louca2018robust,molzahn2018towards}. This may lead to violations of operational constraints and possibly cause voltage collapse. For example, unexpected DC fast charging events or loss of renewable generation can make the system unable to accomplish the load-supporting task, in which case the power flow equations cease to admit a solution~%
    \cite{simpson2016voltage,cui2018voltage}.
    
    Second, previous results do not consider stability issues of DN-OPF solutions. A DC network has rich dynamics contributed by electrical circuits and control systems~%
    \cite{dorfler2018electrical}, which are subject to notable risks of instability associated with the choice of operating point~%
    \cite{7328761,7182770,6031929,zonetti2016tool,7798761}. Many loads in a DC network can be considered as CPLs that are known to have harmful negative impedance effects. If the operating point is not carefully selected, the system can be under-damped or even become unstable~%
    \cite{7328761}. It is worth mentioning that the joint problem of existence and stability of an equilibrium has been studied in control-theoretic work such as~\cite{ikeda1991parametric,wada1998parametric}. However, these prior works primarily focus on developing feasibility and stability conditions. The literature still lacks a control synthesis approach that not only ensures existence and stability of the equilibrium but also designs the location of the equilibrium state, for example, to guarantee the satisfaction of various engineering constraints and achieve economic operations for a DC network.

We propose a stability-constrained robust DN-OPF algorithm to address these limitations. Following DC network operation practices, we focus on a DC network with nonlinear CPLs and controllable voltage sources. We seek to minimize system operational costs by computing setpoints for the sources which rigorously guarantee the following two properties for any loading condition within a specified uncertainty set: 1)~robust feasibility (existence of power flow solutions satisfying operational constraints) and 2) robust stability (local exponential stability of the operating point). 

To provide such guarantees, we formulate a DN-OPF problem that incorporates robust feasibility and stability conditions. Solving this problem is difficult. First, existing stability conditions for DC networks are developed to study given operating points~%
\cite{liu2018robust}; hence, ensuring stability when operating points are decision variables is challenging. Additionally, to ensure robustness, the power flow equations along with the stability conditions need to jointly hold for all uncertainty realizations. This results in a semi-infinite programming (SIP) problem~%
\cite{hettich1993semi} that is generally computationally intractable~%
\cite{mulvey1995robust}.
Tractable reformulations or approximations exist for robust optimization problems when special structures of problem formulation and data uncertainty can be exploited~%
\cite{ben2009robust};
however, there are no standard approaches to deal with the nonconvexity associated with the power flow equations.

\begin{figure}[ht]
	\centering
	\includegraphics[width=0.45\textwidth]{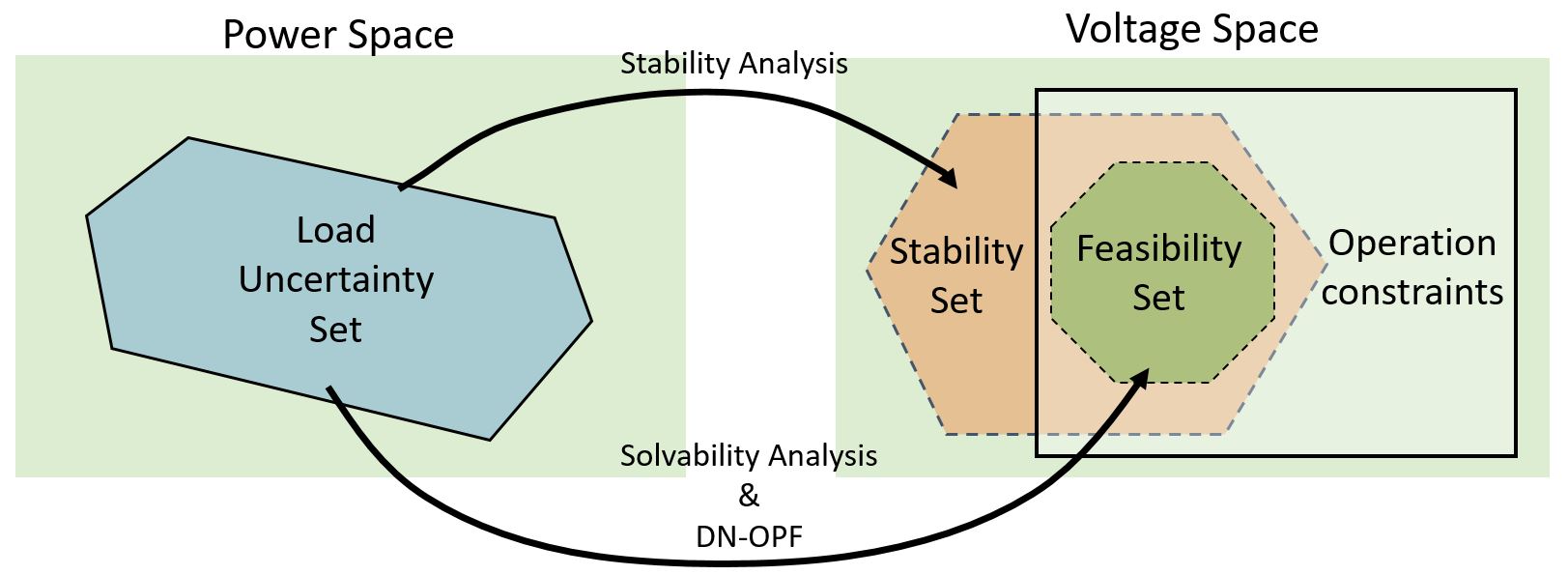}
	\caption{\label{fig:idea}Illustration of the proposed work.}
\end{figure}

The proposed algorithm converts the SIP problem into a tractable formulation that resembles a well-studied DN-OPF problem. The main idea of the proposed work is illustrated in Fig.~\ref{fig:idea}, and we summarize the main technical tasks as follows:
\begin{itemize}
	\item[(1)] We first derive a stability set in the voltage space such that any operating point therein is guaranteed to be stable.
	\item[(2)] We then develop conditions that guarantee the existence of a power flow solution in a feasibility set for any loading condition. We characterize the boundaries of the feasibility set defined by these conditions.
	\item[(3)] Finally, we formulate and solve a tractable problem reminiscent of a \mbox{DN-OPF} problem to ensure that the entirety of the set of operating points lies in the intersection of the stability set and the operational constraints.
\end{itemize} 

The contributions of the paper are summarized as follows:
\begin{itemize}
\item We develop a novel algorithm to reformulate and solve a class of otherwise intractable DN-OPF problems using a synthesis of new DC network stability analysis and power flow feasibility results.
\item We provide a new condition regarding the solvability of DC network power flow equations and new methods to computing a stability set. The solvability condition establishes an explicit bound on power flow solution (bus voltages) as a function of load powers.
\item We provide insights into DC network operations. For example, we provide a rigorous argument for monotonicity in power flow solutions such that the reduction of load power at any bus leads to strictly higher load bus voltages for all load buses.
\end{itemize}

The rest of the paper is organized as follows: First, Section~\ref{sec:prob} introduces the system model and the main problem considered in this paper. Next, Section~\ref{sec:technical} shows the main results of the paper, i.e., a solution algorithm for a robust \mbox{DN-OPF} problem with feasibility and stability guarantees. Section~\ref{sec:simu} then demonstrates the efficacy of the proposed work using simulation case studies. Finally, Section~\ref{sec:sum} concludes the paper and discusses future research directions.

\section{System Modeling and Problem Statement}\label{sec:prob}
\subsection{Notation}
In this paper, we use $\mathbbold{1}$ and $\mathbbold{0}$ to represent vectors of all 1's and 0's of appropriate sizes. Recall that a square matrix $A$ is Hurwitz if all real parts of its eigenvalues are negative. For a vector $v$, let $v_k$ represent its $k$-th element. Let the operator $\diag\{v\}$ yield a diagonal matrix with the vector's components being the diagonal entries. For a real square matrix $A$, $A^{-1}$ denotes its inverse, $A \succ 0$ (resp., $A \succeq 0$) means it is symmetric positive definite (resp., semidefinite), and $A \prec 0$ (resp., $A \preceq 0$) means $-A \succ 0$ (resp., $-A \succeq 0$).

\subsection{DC Power Systems}\label{sec:prob_model}
\begin{figure}[t]
	\centering
	\includegraphics[width=0.35\textwidth]{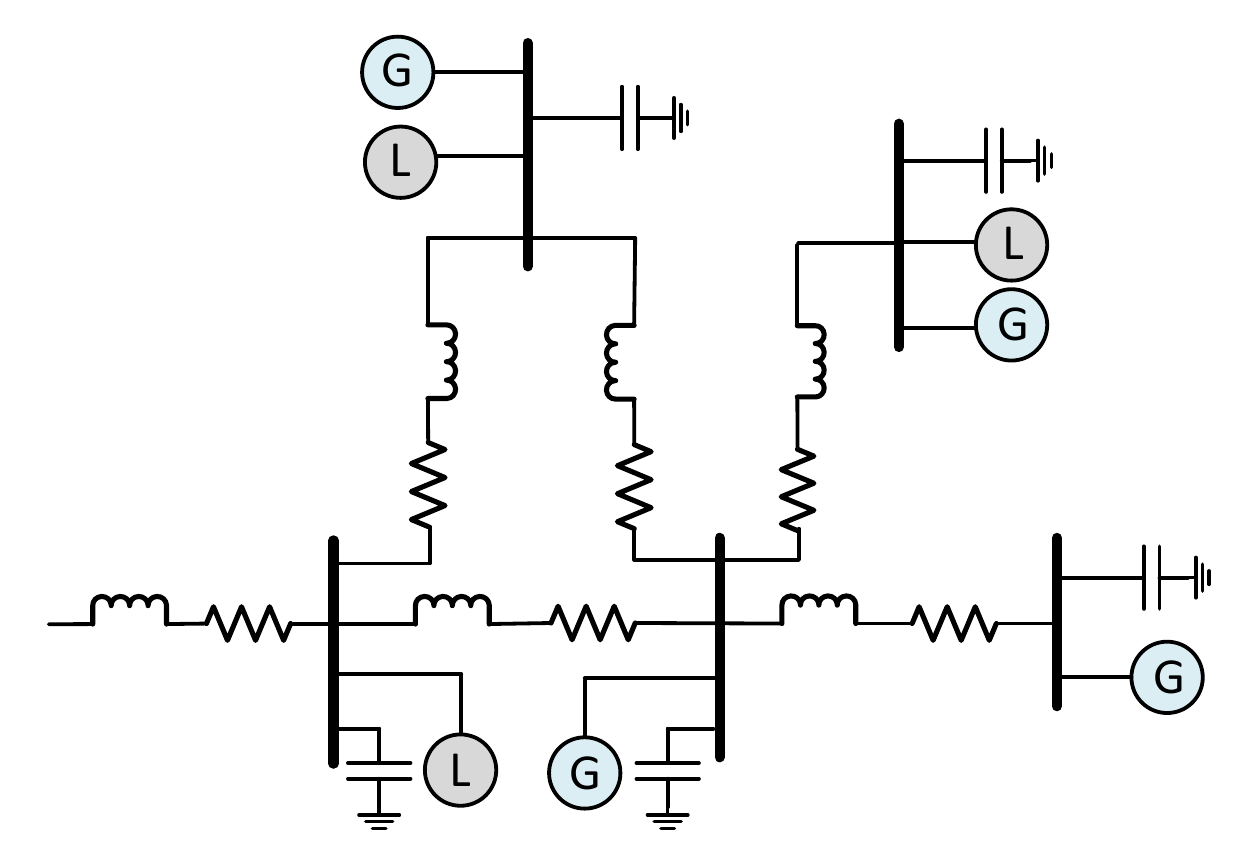}
	\caption{\label{fig:topo}Example DC power network.}
\end{figure}

In this paper, we focus on a DC network with $n_s$ generators, $n_\ell$ loads, and $n_{\mathrm{c}}$ power lines. The total number of these components is $n=n_s+n_{\mathrm{c}}+n_\ell$. Let the index sets of generators, loads, and power lines be $\mathcal{N}_s$, $\mathcal{N}_\ell$, and $\mathcal{E}_{\mathrm{c}}$, respectively. Fig.~\ref{fig:topo} shows an example DC network consisting of lumped $\pi$-equivalent models~%
\cite{7798761} where generators and loads are interconnected via equivalent RLC circuits~%
\cite{dorfler2018electrical}.
\subsubsection{Load and Generator Models}
\begin{figure}[t]
	\centering
	\includegraphics[width=0.4\textwidth]{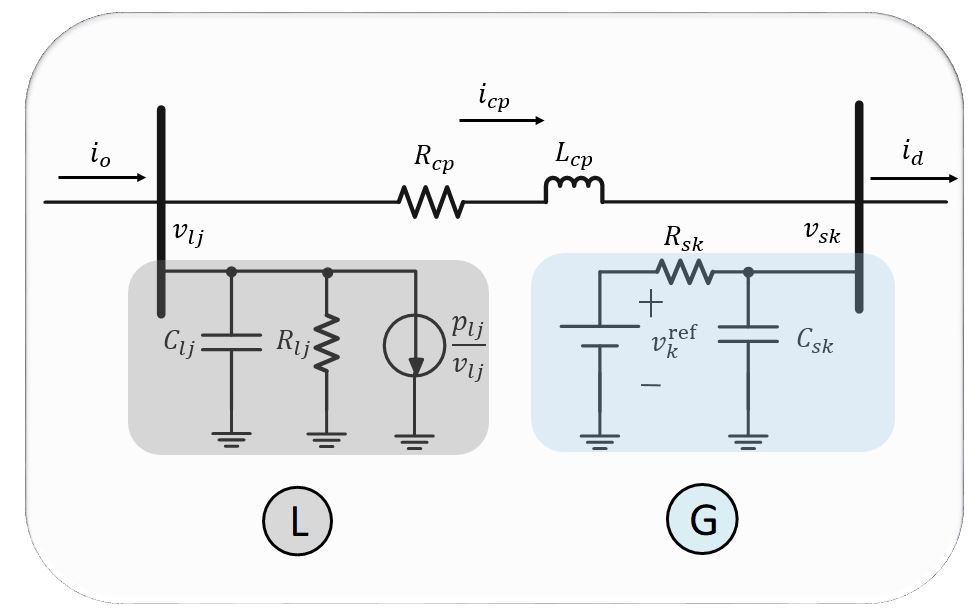}
	\caption{\label{fig:indi_gl}Zoomed-in image of the dynamic circuit.}
\end{figure}

Fig.~\ref{fig:indi_gl} shows a zoomed-in image of one part of the circuit. The lines are represented using a $\pi$-equivalent model~%
\cite{7798761} with a series line resistance and line inductance connected with shunt capacitors at both ends,\footnote{Note that the line resistance, inductance, capacitance are evenly distributed for an ideal power line. The $\pi$-equivalent model is a simplified power line model commonly used for DC network stability analysis~\cite{karlsson2003dc,beerten2013modeling,dorfler2018electrical}.} and each DC bus has a DC bus capacitor and a shunt resistor. The shunt capacitors are in parallel with the DC bus capacitor. For simplicity, we use one composite shunt capacitor to model their joint effects. \color{black}{Suppose the circuit has the $k$-th generator, $p$-th power line, and $j$-th load. Let $i_{o}(t)$ and $i_{d}(t)$ represent the current flowing into and out of the circuit}, and let $i_{\mathrm{c}p}$ represent the current flow in power line $p$. 

Loads are modeled as constant power injections into the network. We let both constant power loads (CPLs) and constant power sources (CPSs) be connected in parallel. These are non-ideal components with parasitic resistances that are represented as lumped resistors. It is well known that CPLs are nonlinear loads and that their associated negative impedance effects are major sources of instabilities in DC networks~\cite{7182770,6415284}. 

For the $j$-th load, let $p_{\ell j}$ represent its power injection into the network, and let $v_{\ell j}$ represent the terminal voltage. At the nominal condition, $p_{\ell j}=p_{\ell j}^*$, where $p_{\ell j}^*$ is a given constant. Each $p_{\ell j}$ is considered to be a perturbation to $p_{\ell j}^*$ that is unknown and bounded within a given uncertainty interval $[\ubar{p}_{\ell j},\bar{p}_{\ell j}]$. Throughout this paper, we let positive $p_{\ell j}$ represent positive power injection (i.e., power generation) at bus $j$. We also let $\bar{p}_{\ell j} \geq 0$ to allow loads to be turned off completely. In this paper, the load profile $p_\ell$ is modeled as an uncertain vector bounded by element-wise interval constraints. Let $\mathcal{P}_\ell$ be the interval uncertainty set, that is, $\mathcal{P}_\ell = \{p_\ell: p_{\ell j}\in [\ubar{p}_{\ell j},\bar{p}_{\ell j}], j \in \mathcal{N}_\ell \}$. This modeling choice is appropriate for typical DC networks in conventional settings. DC networks are often used for medium- and low-voltage power distribution systems~\cite{dragivcevic2016dc}. For such systems, the passive load variability is generally spatially uncorrelated. We also note that the time-scale of the RLC dynamics of concern is usually milliseconds to sub-seconds, while the time-scale of load variations usually range from seconds to minutes. Hence, the load power profile is considered constant in the dynamical model with uncorrelated interval uncertainty.

\color{black}
For the $k$-th source, let $V^{\mathrm{ref}}_k$ be the controllable voltage setpoint, $v_{sk}$ be the external generator voltage, and $R_{sk}$, $C_{sk}$ represent the source resistance and capacitance, respectively. We impose operational constraints on controllable voltage setpoints such that vector $V^{\mathrm{ref}}$, which includes all voltage setpoints, needs to lie within a given convex constraint set~$\mathcal{V}^{\mathrm{ref}}$. 
\begin{rem}[Generator Model]
	Generators are modeled as non-ideal voltage sources~%
	\cite{7798761} that are in series with equivalent resistors. The voltage set-points can be controlled by a grid operator. We assume that proper low-level controllers~%
	\cite{dragivcevic2016dc} have been employed to regulate the terminal voltage of a generator to track a reference setpoint. Consequently, the generator can automatically vary power outputs to respond to changing loading conditions. The generator internal dynamics, including those from low-level controllers and electromechanical transients, are not considered in this paper, and we mainly focus on the network dynamics contributed by electromagnetic transients in the stability analysis. Nevertheless, the main results of the paper can be extended to various generator dynamics as well. It is worth mentioning that these generators are different from the CPSs. The power output of a CPS is uncertain, and its terminal voltage cannot be directly regulated. For example, solar converters working in Maximum Power Point Tracking mode can be considered as a CPS. The converter is in power-control mode to regulate the power output with respect to uncertain solar irradiance~\cite{ram2017comprehensive}.
\end{rem} 	
\begin{rem}[Extension]
	The main results of this paper can be extended to DC networks with other generator and load models. For example, constant-current and constant-impedance loads are linear and can be easily incorporated in the model. Additionally, for generators with V-I droop control~%
	\cite{dragivcevic2016dc}, the voltage setpoint can be considered as the droop reference and the droop gains can be equivalently modeled as virtual impedances, which only changes the parameters of the RLC circuits.
\end{rem} 
\begin{rem}[Implication of Uncertainty Model]
    We elaborate the implication of the uncertainty model on optimization formulation and computation in the following. Modeling the uncertainty as a vector with interval bounds is the only known scenario under which the resulting uncertain SDP constraint, which is used to ensure robust stability, admits a tractable approximation with guaranteed tightness factor~\cite{ben2009robust}. It is worth mentioning that more general uncertainty sets can be directly adopted in the proposed framework as long as a tractable approximation of the uncertain SDP constraint exists, regardless of explicit knowledge of approximation quality. For instance, \cite{zhen2017robust} shows that tractable approximations of uncertain SDP problems can be derived for any polyhedral uncertainty set. It then follows that the proposed approach can be applied to problems with polyhedral uncertainty sets. Since ellipsoids can be well approximated by polytopes, ellipsoidal uncertainty sets, which are commonly employed for the modeling of uncertainty set of renewable power profile, can be handled in the proposed framework as well. However, we use the interval uncertainty model in this paper due to its relevance for DC networks and brevity of exposition.
\end{rem}

\subsubsection{Dynamic Network Model} 
Sources and loads are connected to DC buses. The buses form a connected graph where a bus is a node and an edge is a $\pi$-equivalent power line. 

We exemplify the modeling approach using the circuit shown in Fig.~%
\ref{fig:indi_gl}. The state variables of the example circuit are the voltages of the capacitors and the currents through the inductors, namely, $v_{sk}(t)$, $v_{\ell j}(t)$, and $i_{\mathrm{c}p}(t)$. The design variables are the output voltages of the sources, $V^{\mathrm{ref}}_k$, $\forall k\in\mathcal{N}_s$. The dynamics of the circuit are represented by the following model using Kirchhoff's current and voltage laws, 

\begin{subequations} \label{eq:mdl_circ}
	\begin{align}
	\label{eq:mdl_circ_a}
	\frac{di_{\mathrm{c}p}(t)}{dt} &= \frac{1}{L_{\mathrm{c}p}} \left( v_{\ell j}(t)-R_{\mathrm{c}p}i_{\mathrm{c}p}(t)-v_{sk}(t) \right), \\
	\label{eq:mdl_circ_b} 
	\frac{dv_{sk}(t)}{dt} &= \frac{1}{C_{sk}} \left( \frac{V^{\mathrm{ref}}_k-v_{sk}(t)}{R_{sk}} - i_{d}(t) + i_{\mathrm{c}p}(t) \right), \\
	\label{eq:mdl_circ_c}
	\frac{dv_{\ell j}(t)}{dt}\! &= \! \frac{1}{C_{\ell j}} \left(-\frac{v_{\ell j}(t)}{R_{\ell j}} + i_{o}(t) - i_{\mathrm{c}p}(t) + \frac{p_{\ell j}}{v_{\ell j}(t)} \right).
	\end{align}
\end{subequations}
Equations~\eqref{eq:mdl_circ_a} and~\eqref{eq:mdl_circ_b} characterize the behavior of the power line and the source. They are linear in the state and design variables. However,~\eqref{eq:mdl_circ_c} is nonlinear due to the term, $p_{\ell j}/v_{\ell j}(t)$. Recall that $i_{o}(t)$ and $i_{d}(t)$ represent aggregate currents flowing from or into the rest of the network. Each of them is a linear combination of the line currents injections into the load bus or the generator bus. 

The modeling approach can be applied to the entire system. By dropping the subscripts indicating variable indices, $p_\ell, v_\ell, v_s, i_{\mathrm{c}}$, $V^{\mathrm{ref}}$ represent the vectors of load powers, load voltages, generator external voltages, power line currents, and controllable voltage setpoints, respectively. Let $x=[i_{\mathrm{c}}^\top,v_{s}^\top,v_{l}^\top]^\top$ be the vector of state variables and $h(x,p_\ell)=[p_{\ell1}/v_{\ell1}, \ldots, p_{\ell n_\ell}/v_{\ell n_\ell}]^\top$.

With the above description and notation, the overall dynamics of the DC grid can be written as follows:
\begin{equation}
    \dot{x}(t)=Ax(t)+BV^{\mathrm{ref}}+Ch(x(t),p_\ell), \quad p_\ell\in \mathcal{P}_\ell, \label{eq:gencls}
\end{equation}
where $A\in \mathbb{R}^{n\times n}$, $B\in \mathbb{R}^{n\times n_s}$, and $C\in \mathbb{R}^{n\times n_\ell}$ are constant matrices that are determined by the network topology and RLC circuit parameters through similar methods to those in~%
\cite{liu2018robust}. This is a well-accepted model for DC network stability studies and has been applied to a variety of applications~\cite{dragivcevic2016dc,kalcon2012dctrans,1658410} (e.g., analyses of DC transmission system dynamics~%
\cite{kalcon2012dctrans}).

\begin{rem}[Control Dynamics]
	So far, we have discussed a general model consisting of additive linear and nonlinear parts to represent the circuit dynamics. The model can be extended to represent control dynamics in a DC network as well. For example, averaging proportional-integral (DAPI) algorithms~\cite{molzahn2017survey,7061540} have recently been developed to control a DC network. They introduce new linear dynamics to our model and the developed results can be extended to incorporate them as well.
\end{rem}

Let $x^\mathrm{e}(p_\ell,V^\mathrm{ref}) = [(I_{\mathrm{c}})^\top, (V_s)^\top, (V_\ell)^\top]^\top \in \mathbb{R}^{n}$ be an equilibrium of~\eqref{eq:gencls}. Notice that the equilibrium is a function of the load power profile and source voltage setpoint. For notational simplicity, we drop the arguments in $I_{\mathrm{c}}$, $V_s$, and $V_\ell$ here and denote the equilibrium by $x^\mathrm{e}$ in the rest of the paper.

When $V_{\ell k} \neq 0, \forall k \in \mathcal{N}_\ell$, the linearized Jacobian matrix of system~\eqref{eq:gencls} with respect to $x^\mathrm{e}$ is given as follows:

\begin{equation}
    J(V_\ell ,p_\ell) = A - \sideset{}{_{k \in \mathcal{N}_\ell}}\sum\frac{p_{\ell k}}{V^2_{\ell k}}e_k^{\vphantom{\top}} e_k^\top,\label{eq:jaco}
\end{equation}
where $e_k \in \mathbb{R}^{n}$ is a basis vector with the \mbox{$(n_{\mathrm{c}}+n_s+k)$-th} element being $1/\sqrt{C_{\ell k}}$. The Jacobian matrix is an affine function in each term $p_{\ell k}/V^2_{\ell k}$. The form of this matrix shows that the local stability of an operating point depends on both the CPL power and the steady-state CPL voltage. 

From basic control theory \cite[Thm. 4.6]{khalil2002nonlinear}, an equilibrium $x^\mathrm{e}$ of \eqref{eq:gencls} is locally exponentially stable if there exists a real $n \times n$ positive definite matrix $P$ that satisfies the following condition:
\begin{equation}
PJ(V_\ell, p_\ell)+J(V_\ell, p_\ell)^\top P \prec 0. \label{eq:con_stab}
\end{equation}
If $p_\ell$ and $V_\ell$ are given, this condition is a linear matrix inequality (LMI) constraint. However, in our problem, $p_\ell$ is uncertain, $V_\ell$ is a variable to be determined, and the coupling between $V_\ell$ and $P$ is non-polynomial.

\subsubsection{Power Flow Model} The power flow model describes the steady-state behavior at an operating point of a DC network. The power flow model is obtained by setting the left-hand side of~\eqref{eq:gencls} to $\mathbbold{0}$ and rearranging terms as
\begin{align}
p_\ell=\diag\{V_\ell\}\left( Y_{\ell \ell}V_\ell+Y_{\ell s}V^{\mathrm{ref}}\right), \label{eq:pfe_ori}
\end{align}
where the connectivity between CPL--source and CPL--CPL are described by two admittance matrices $Y_{\ell s}\in \mathbb{R}^{n_\ell\times n_s}$ and $Y_{\ell\ell}\in \mathbb{R}^{n_\ell\times n_\ell}$~%
\cite{liu2018existence}, which are submatrices of the system admittance matrix $Y$. Equation~\eqref{eq:pfe_ori} is quadratic in state variables $V_\ell$ and bilinear in design variables $V^{\mathrm{ref}}$ and state variables $V_\ell$. In general, the power flow model~\eqref{eq:pfe_ori} usually introduces computational challenges owing to its nonconvexity~\cite{gan2014optimal}.

In addition, the system at steady state needs to satisfy operational constraints. In this paper, we require that $V_\ell \in \mathcal{V}_\ell$ and $V^{\mathrm{ref}} \in \mathcal{V}^{\mathrm{ref}}$. Both sets are convex sets that represent system operational requirements such as upper and lower voltage bounds. We limit our presentation to only consider the constraints related to the load voltages and generator setpoints, which are directly relevant to the system stability, in order to simplify the paper's discussion. Other variables like the currents are linear functions of the load voltages and generator voltage setpoints. The proposed algorithm can be easily extended to incorporate constraints on these variables.

\subsection{Problem Statement}\label{sec:prob_problem}
From the models discussed above, a poorly designed $V^{\mathrm{ref}}$ may 1)~result in violations of operational constraints; 2) cause local instability for operating points; and 3)~lead to infeasibility of \eqref{eq:pfe_ori} or even loss of equilibrium altogether. 

The goal of this work is to choose the value of $V^{\mathrm{ref}}$ which minimizes operating costs while guaranteeing that the system is robustly feasible and stable. We make the terms \emph{robustly feasible} and \emph{robustly stable} precise in Definition~\ref{def:feas+stab} below:

\begin{definition} \label{def:feas+stab}
	Given a generator voltage setpoint $V^{\mathrm{ref}}$, system~\eqref{eq:gencls} is said to be \emph{robustly feasible} if, for every $p_\ell \in \mathcal{P}_\ell$, the system admits an equilibrium $x^\mathrm{e}$ which satisfies all operational constraints. The system is said to be \textit{robustly stable} if, for every $p_\ell \in \mathcal{P}_\ell$, there exists a corresponding $V_\ell$ such that the Jacobian $J(V_\ell ,p_\ell)$ is Hurwitz.
\end{definition}

\begin{rem}\label{rem:feas+stab}
   Although robust control problems have been extensively studied in the literature~\cite{leitmann1979guaranteed,doyle1989state}, a special feature of our problem is that the existence, location, and local stability of the equilibrium state is dependent on the interplay between the parameter uncertainty and the control input, as opposed to a common assumption that the equilibrium is fixed at the origin~\cite{ikeda1991parametric}.
\end{rem}

Desirable operating points for power systems are usually computed by solving OPF problems~%
\cite{cain2012history}. Recently, OPF problems for DC networks (DN-OPF) have been a particular research focus~%
\cite{gan2014optimal,li2018optimal,montoya2018linear}. The formulation of existing \mbox{DN-OPF} problems can be summarized as follows:
\begin{subequations}\label{eq:prob_ori}
	\begin{align}
	\textbf{DN-OPF$^*$: }\min_{\substack{V^{\mathrm{ref}}\in \mathcal{V}^{\mathrm{ref}}}}  & \, f(V^{\mathrm{ref}},V_{\ell}^*), \\
	\text{s.t.} \quad & p^*_\ell \!=\!\diag\{V_\ell^*\}\!\left( Y_{\ell \ell}V^*_\ell \! + \! Y_{\ell s}V^{\mathrm{ref}}\right)\!,\label{eq:con_pfe_ori} \\  
	& V^*_\ell \in \mathcal{V}_\ell, \label{eq:con_pfe_ori_2}
	\end{align}%
\end{subequations}%
where $p^*_\ell\in \mathbb{R}^{n_\ell}$ is the nominal CPL power profile, $V_\ell^*$ is the steady-state load voltage at the nominal load condition, and \mbox{$f:\mathbb{R}^{n_s}\times \mathbb{R}^{n_\ell}\to \mathbb{R}$} is a possibly nonconvex cost function that usually represents the operating cost (e.g., power loss or generation cost). 

Note that in classical OPF formulations, the cost functions typically consider the generation costs. Since our main results do not depend on the structure of the cost function, we do not explicitly specify the cost functions in this paper in order to simplify our discussion. Our approach can easily accommodate typical OPF cost functions. For example, the case studies in Section~\ref{sec:simu} minimize generation costs. The total generation cost is formulated by summing the products of each generator's cost coefficient with its power output. Since the generator power outputs can be found using power flow equations similar to~\eqref{eq:pfe_ori}, this choice of objective $f$ can be written in terms of the voltages.

Recently, effective methods have been developed to solve the DN-OPF problem~\eqref{eq:prob_ori} using approximation and convex relaxation techniques~%
\cite{gan2014optimal,li2018optimal,montoya2018linear}. However, problem~\eqref{eq:prob_ori} only considers a fixed loading condition and does not explicitly consider system stability. If the actual load is different from the nominal load, the system's operating point may be unexpected and possibly even unstable.

To address these limitations, we focus on the following problem with explicit constraints guaranteeing robust feasibility and robust stability:
\begin{subequations}\label{eq:prob1}
	\begin{align}
	\noalign{\textbf{R. DN-OPF SIP: }}
	\min_{\substack{V^{\mathrm{ref}}\in \mathcal{V}^{\mathrm{ref}} \\ P \succeq 0}} \quad & f(V^{\mathrm{ref}},V^*_\ell), \\
	\text{s.t.} \quad & (\forall p_\ell\in \mathcal{P}_\ell) \nonumber \\
	& P J(V_\ell(p_\ell), p_\ell) + J(V_\ell(p_\ell), p_\ell)^\top P \prec 0, \label{eq:prob1:b} \\
	& p_\ell=\diag\{V_\ell(p_\ell)\}\left( Y_{\ell \ell}V_\ell(p_\ell) + Y_{\ell s}V^{\mathrm{ref}} \right), \label{eq:prob1:c} \\
	& V_\ell(p_\ell) \in \mathcal{V}_\ell,\quad~\eqref{eq:con_pfe_ori}, \quad ~\eqref{eq:con_pfe_ori_2}. \label{eq:con_seminf}
	\end{align} 
\end{subequations}
Compared to problem~\eqref{eq:prob_ori}, we add the sufficient stability condition from~\eqref{eq:con_stab} in order to ensure robust stability. We also require all constraints to hold for all $p_\ell\in \mathcal{P}_\ell$ in order to ensure robust feasibility in the presence of uncertainty. Note that $V_{\ell}(p_\ell)$ is an implicit variable; with given $V^{\mathrm{ref}}$, $V_{\ell}(p_\ell)$ is a function of $p_\ell$.


Problem \eqref{eq:prob1} is a robust optimization problem involving an LMI constraint with structured uncertainty as in~\eqref{eq:prob1:b} and nonconvex constraints as in~\eqref{eq:prob1:c}. While tractable reformulations or safe approximations have been identified for some robust conic programming problems with special uncertainty sets~\cite{ben2009robust}, there are no known general tractable reformulations or safe approximations for robust nonconvex optimization problems in the form of \eqref{eq:prob1}. In fact, it has been shown in~\cite{lehmann2016ac} that the AC-OPF problem is NP-hard even in the deterministic case. Although DN-OPF is simpler than its AC counterpart, there is no efficient solver with global optimality guarantee as far as we know. Even though standard OPF problems can often be solved with high quality by modern interior point-based solvers despite their theoretical computational complexity, the lack of tractable reformulations and approximations makes solving its robust counterpart difficult. As a result, existing robust OPF formulations almost exclusively adopt power flow models that employ convex relaxation techniques instead of the original nonconvex formulation~\cite{lorca2018adaptive, louca2018robust, yang2020robust}. To find a tractable way to search for a feasible solution to~\eqref{eq:prob1}, we derive an efficient convex inner approximation of the feasible region of problem~\eqref{eq:prob1} in the next section.

\section{Tractable DN-OPF with Robust Feasibility and Stability Guarantees} \label{sec:technical}

In this section, we derive a computationally tractable optimization problem whose feasible region is a convex inner approximation to that of the original problem~\eqref{eq:prob1}. As illustrated in Fig.~\ref{fig:idea}, the development of the inner approximation consists of three main steps. First, we solve a series of tractable SDPs to find a stability set in $\mathcal{V}_\ell$ such that for any load in $\mathcal{P}_\ell$, the corresponding equilibrium is locally exponentially stable if its load voltages lie in the stability set. Second, based on a power flow solvability technique, we derive a sufficient condition on voltage setpoints $V^\mathrm{ref}$ such that a feasible power flow solution within operational constraints and stability set exists for every $p_\ell \in \mathcal{P}_\ell$. Third, we formulate a DN-OPF problem which optimizes $V^{\mathrm{ref}}$ over $\mathcal{V}^{\mathrm{ref}}$ while satisfying the sufficient condition derived in the second step.
\subsection{Robust Stability Set}\label{sec:robstab}
The stability set is the feasibility region of the stability condition~\eqref{eq:con_stab}. Due to the infinite number of constraints and the non-polynomial structure, this region is difficult to characterize. Motivated by power system operational constraints, this section describes an interval set which inner approximates this region. 

Let $\ubar{V}^s_{\ell j}$ and $\bar{V}^s_{\ell j} \in \mathbb{R}^{n_\ell}_+$ represent the lower and upper bounds of an interval set of $V_\ell$, denoted as $\mathcal{V}^s_\ell=\left\{ V_\ell :\; \ubar{V}^s_{\ell j} \leq V_{\ell j} \leq \bar{V}^s_{\ell j}, \forall j \in \mathcal{N}_\ell \right\}$. We term $\mathcal{V}^s_\ell$ a ``robust stability set'' when the following definition applies:
\begin{definition} \label{def:robstab}
	A set $\mathcal{V}^s_\ell$ is called a \emph{robust stability set} if there exists a positive definite matrix $P$ such that the following inequality is satisfied for all $V_\ell \in \mathcal{V}^s_\ell$ and all $p_\ell \in \mathcal{P}_\ell$:
	\begin{equation}\label{eq:prob_nonpoly}
	 PJ\left(V_\ell,p_\ell\right) + J\left( V_\ell,p_\ell\right)^\top P \prec 0.\vspace*{0.5em}
	\end{equation} 
\end{definition}

From this definition, if an operating point lies in a robust stability set under load uncertainty set $\mathcal{P}_\ell$, it is locally exponentially stable regardless of specific realization $p_\ell \in \mathcal{P}_\ell$. As shown below, this provides us with the flexibility to remove the coupling between equilibrium and load profiles.

The interval robust stability set facilitates efficient optimization formulation to solve problem~\eqref{eq:prob1}. In the following subsections, we develop rigorous upper and lower bounds on the load bus voltages under any load uncertainty realization given the generator voltage setpoints $V^\mathrm{ref}$. To certify robust stability, one only needs to check that the load bus interval is included in the interval robust stability set, which can be performed in a computationally efficient manner.

\subsubsection{Interpolation}
The matrix $J(V_\ell,p_\ell)$ has the following two special structures: first, the variable $V_{\ell k}$ and the uncertain parameter $p_{\ell k}$ only exist in pairs on the diagonal entries in the form $-p_{\ell k}/V_{\ell k}^2$; second, each composite term $-p_{\ell k}/V_{\ell k}^2$ only appears once in the matrix. This provides the possibility of applying an interpolation method to replace each $-p_{\ell k}/V_{\ell k}^2$ with a new variable.

Let $\delta_k = - p_{\ell k}/V_{\ell k}^2$ and $\delta = [\delta_1,\cdots,\delta_{n_\ell}]^\top$. When $p_{\ell k}$ and $V_{\ell k}$ are subject to box constraints, the vector $\delta$ is contained in an interval set as well. Let $\Delta \triangleq \{\delta:\ubar{\delta}_k\leq \delta_k\leq \bar{\delta}_k\}$. We call $\delta$ the diagonal perturbation to the system Jacobian and call $\Delta$ the diagonal perturbation set. Since $p_{\ell k}$ can be positive or negative, we let $\ubar{\delta}_k < 0$ and $\bar{\delta}_k > 0$.

With the above discussed definition, substituting~$\delta_k = -p_{\ell k}/V^2_{\ell k}$ into~\eqref{eq:jaco} yields a new expression for the system Jacobian, $J(\delta)=A+\sum_{k \in \mathcal{N}_\ell}\delta_k^{\vphantom{\top}} e_k^{\vphantom{\top}} e^\top_k$. Hence, the linearized system matrix is now subject to an affine interval parameter uncertainty. In the following, we denote this matrix as $J(\delta)$. We are interested in finding a diagonal perturbation set whose every element makes the matrix Hurwitz stable, i.e., satisfy the following inequalities:
\begin{equation}
    PJ(\delta)+J^\top(\delta)P \preceq 0, \quad \forall \delta \in \Delta.\label{eq:stab_ori}
\end{equation}
Notice that once this $\Delta$ can be found, we can translate it into the desired robust stability set.

For a given $\Delta$, there exist multiple methods to certify whether~\eqref{eq:stab_ori} is satisfied~\cite{dorfler2018electrical,liu2018robust}. Most existing methods check multiple ``critical scenarios'' to certify constraint satisfaction for all scenarios. The existing methods have issues with conservativeness or computational tractability. For example, the condition in~\cite{dorfler2018electrical} tests whether all diagonal elements of the Jacobian matrix are negative, which cannot be satisfied in our case when $\delta_k$ is positive. A sufficient condition based on LMI feasibility testing is proposed in~\cite{liu2018robust} that involves $2^{n_\ell}$ LMI constraints. While numerical tests reveal that the condition in~\cite{liu2018robust} has advantages with respect to limited conservativeness, practical applicability of this condition is challenging since the number of LMI constraints is exponentially dependent on the dimension of the uncertainty. 

Since~\eqref{eq:stab_ori} is a \emph{robust semidefinite programming} problem with interval uncertainty set,\footnote{In a robust semidefinite programming problem, the constraint is bilinear in the uncertainty and the decision variable while it is an LMI in the decision variable if the uncertainties are known constants.} Theorem~9.1.2 of~\cite{ben2009robust} can be applied to develop a new condition for DC network stability analysis:
\begin{lem}~\label{lem:LMIs}
Given $\Delta$, if there exists $P \succ 0$, $N \prec 0$, and $n_\ell$ positive scalars $\lambda_1$, $\cdots$, $\lambda_{n_\ell}$ that satisfy the following LMI conditions, $J(\delta)$ is always Hurwitz stable for all $\delta \in \Delta$:
\begin{subequations}\label{eq:LMItest_poly}
	\begin{align}
	&\left[
	    \begin{array}{cccc}
	        N\!+\!\sum_{k=1}^{n_\ell}\lambda_k e_k^{\vphantom{\top}} e_k^\top \left(\frac{\bar{\delta}_k-\ubar{\delta}_k}{2}\right)^2 
	        & P e_1 & \cdots & P e_{n_\ell} \\
	        e_1^\top P & -\lambda_1\\
	        \vdots & & \ddots\\
	        e_{n_\ell}^\top P & & & -\lambda_{n_\ell}
	    \end{array}
	\right] \preceq 0, \label{eq:LMItest_poly_1}\\
	&N \succeq P (A + \sum_{k=1}^{n_\ell} e_k^{\vphantom{\top}} e_k^\top \frac{\bar{\delta}_k + \ubar{\delta}_k}{2} ) + (A + \sum_{k=1}^{n_\ell} e_k^{\vphantom{\top}} e_k^\top \frac{\bar{\delta}_k + \ubar{\delta}_k}{2} )^\top P.\label{eq:LMItest_poly_2}
	\end{align}
\end{subequations}    
\end{lem}
The condition in Lemma~\ref{lem:LMIs} only involves two LMIs. The decision variables are two $n \times n$ matrices, $P$ and $N$, as well as $n_\ell$ scalars $\lambda_j$. In total, the condition has $2n^2+n_\ell$ free scalar variables. Since this number is polynomially dependent on $n$ and $n_\ell$, the condition has reasonable scalability. 

The reduction in computational complexity may induce concerns regarding conservativeness. One method to evaluate the conservativeness is to compare the volume of the largest sets that the conditions can certify. Compared to Lemma~1 of~\cite{liu2018robust} where the number of LMIs is exponential in the number of loads, numerical tests show that the proposed conditions can certify a set with a volume above 95$\%$ of the largest set certifiable by exponentially many constraints. For example, we apply our results to the example DC microgrid detailed in Case Study 2 of~\cite{liu2018robust}. The simulation results are shown in Table~\ref{table:csv_1}. Notice that in the case study, we consider all pure-load buses and fix the voltage lower bounds, hence the volume of the stability set is indexed by the load power. Higher loads correspond to larger sets and a reduction in the condition's conservativeness.
%
%
Note that the proposed results can certify a set with a volume over $99.5\%$ relative to the condition with exponentially many LMIs in~\cite{liu2018robust} and shows significant improvements compared to the condition with polynominally many LMIs in Proposition 1 of~\cite{liu2018robust}.

\begin{table}[hb]
\caption{Volume of Largest Certifiable Robust Stability Set}
\label{table:csv_1}
\begin{center}
\begin{tabular}{cccc}
\hline
& Lemma~\ref{lem:LMIs} &  Expo. LMIs~\cite{liu2018robust} & Poly. LMIs~\cite{liu2018robust} \\
\hline
Highest Load (kW) & $-19.87$ & $-19.96$ & $-18.23$ \\ 
\hline
\end{tabular}
\end{center}
\end{table}

\subsubsection{Computing the Robust Stability Set}
With the proposed stability condition, we are equipped with a tractable method to check whether a given $\Delta$ satisfies~\eqref{eq:stab_ori}. We would like to find a set $\Delta$ with a larger volume, as it can be translated into a robust stability set with a larger volume as well. The volume of $\Delta$ is determined by its vertices, hence increasing the volume essentially involves adjusting these vertices. We use a line search method to accomplish this goal. 

A subsequent question concerns selecting an initial guess for $\Delta$ that reduces computational efforts in finding a larger stability set. Rather than arbitrarily choosing an initial guess, we make one suggestion of the initial guess that may find the largest interval stability set in a single shot. This guess covers the entire domain of $\delta$ with respect to all possible power flow solutions. The details are provided in Appendix~\ref{sect:app:set}.

Let the initial guess set be denoted as $\Delta_0$. The volume of $\Delta$ can be adjusted by introducing a positive scaling factor $\alpha$ to all the vertices. We denote the set after adjustment as $\Delta = \alpha\Delta_0$. We want to find the largest $\alpha$ that satisfies the condition of Lemma~\ref{lem:LMIs}. This value is found by solving the following generalized eigenvalue problem (GEVP):
\begin{subequations}\label{eq:prob_gevp}
	\begin{align}
	&\textbf{GEVP: }\quad \max_{\alpha>0, P\succ 0, N\prec 0,\lambda > 0}{\;\alpha} \quad \text{s. t.}\\
	&\left[
	    \begin{array}{cccc}
	        N\!+\!\alpha^2\sum_{k=1}^{n_\ell}\lambda_k e_k^{\vphantom{\top}} e_k^\top \left(\frac{\bar{\delta}_k-\ubar{\delta}_k}{2}\right)^2 
	        & P e_1 & \cdots & P e_{n_\ell}\! \\
	        e_1^\top P & -\lambda_1\\
	        \vdots & & \ddots\\
	        e_{n_\ell}^\top P & & & -\lambda_{n_\ell}\!
	    \end{array}
	\right] \!\preceq  0, \label{eq:gevp_1}\\
     & N \succeq  P(A + \sum_{k=1}^{n_\ell} e_k^{\vphantom{\top}} e_k^\top \frac{\bar{\delta}_k + \ubar{\delta}_k}{2} ) + (A + \sum_{k=1}^{n_\ell} e_k^{\vphantom{\top}} e_k^\top \frac{\bar{\delta}_k + \ubar{\delta}_k}{2} )^\top  P.
	\end{align}
\end{subequations}  

A solution that is arbitrarily close to the global optimum can be found for the GEVP problem~\eqref{eq:prob_gevp} since it is a quasi-convex problem~\cite{boyd2004convex}. For a solution $\alpha$ of~\eqref{eq:prob_gevp}, we are endowed with a robust stability set as described below.
\begin{prop}~\label{prop:stab}
	Given $\Delta_0$, if $\alpha$ is a solution of~\eqref{eq:prob_gevp}, $\mathcal{V}^{s}_\ell$ defined in the following is a robust stability set:
	\begin{equation}
	\mathcal{V}^{\mathrm{s}}_\ell = \left\{V_\ell:V_{\ell k}\geq V^{\mathrm{s}}_{\ell k}, \forall k \in \mathcal{N}_\ell\right\},\nonumber
	\end{equation}
	where $(V^{\mathrm{s}}_{\ell k})^2 = \max\{-\bar{p}_{\ell k}/(\alpha\ubar{\delta}_k),\,-\ubar{p}_{\ell k}/(\alpha \bar{\delta}_{k})\}$.
\end{prop}
\begin{proof}
     Available in Appendix~\ref{sect:app:prop_stab}
\end{proof}
In Proposition~\ref{prop:stab}, we find the lower bound for the steady-state voltage to ensure robust stability. This bound is in line with engineering observations for DC grid stability: with larger load power, the system should be operated at higher voltage levels to reduce risks of instability.

\begin{rem}~\label{rem:whygood}
    The value of $\alpha$ is critical for determining the robust stability set. As the set includes all relevant uncertainty realizations, we do not need to test if the system is stable for uncertainty outside of the set. Thus, the value of $\alpha$ is upper bounded by 1. When the solution is obtained as 1, the largest interval stability set can be directly found. Otherwise, a line search algorithm can be developed to approximately solve~\eqref{eq:prob_gevp}. In the line search process, a finite sequence of semidefinite programming (SDP) problems need to be solved. There exist numerically efficient algorithms to solve SDP problems~\cite{boyd2004convex}, and thus the line search process is tractable.
\end{rem}

\subsection{Solvability Condition}\label{sec:robfeas}
As shown in Fig.~\ref{fig:idea}, after the robust stability set is found, our next task is to ensure: a) the system operates in this set and b) the system complies with all other operational constraints. This task is equivalent to the following robust feasibility problem: we need to compute a voltage setpoint $V^{\text{ref}} \in \mathcal{V}^{\text{ref}}$ that guarantees the existence of power flow solutions (i.e., solutions to~\eqref{eq:pfe_ori}) that are in $\mathcal{V}^s_\ell \cap \mathcal{V}_\ell$
for all $p_\ell \in \mathcal{P}_\ell$. 


The problem of certifying the existence and characterizing the range of power flow solutions under uncertain power injections has been studied for the AC power flow model through convex restriction~\cite{Lee19}. Due to the intrinsic nonconvexity of AC power flow model, the convexification approach may be conservative. On the other hand, the geometry of DC power flow equations is considerably simpler than its AC counterpart (cf. \cite{Jeeninga2020dc}) so better results may be expected. In fact, we show in this section that for the given interval uncertainty set of power injections, there is an efficient way to \emph{exactly} certify solution existence and characterize their bounds.

We first introduce the following Lemma from \cite{dvijotham2017high}.
\begin{lem}[\hspace{1sp}{\cite[Thm. 3]{dvijotham2017high}}] \label{lem:high}
    Given $p_\ell$ and $V^\mathrm{ref} > \mathbbold{0}$, if \eqref{eq:pfe_ori} is solvable, then there exists a \emph{high-voltage solution} $V_\ell^* > \mathbbold{0}$ to \eqref{eq:pfe_ori} such that $V_\ell^* \ge V_\ell$ for all other solutions $V_\ell$ to \eqref{eq:pfe_ori}. 
\end{lem}

The next Lemma shows the existence and bounds of power flow solutions for any loading condition $p_\ell \in \mathcal{P}_\ell$ when the \emph{high-voltage solutions} for the extreme loading conditions are known.
\begin{lem} \label{lem:exist}
    Let $V_{\ell -}$ be the high-voltage solution for $p_\ell = \ubar{p}_\ell$, then \eqref{eq:pfe_ori} is solvable for all $p_\ell \in \mathcal{P}_\ell$ and the high-voltage solution for any $p_\ell \in \mathcal{P}_\ell$ satisfies $V_\ell(p_\ell) \in [V_{\ell -}, V_{\ell +}]$ where $V_{\ell +}$ is the high-voltage solution for $p_\ell = \bar{p}_\ell$.
\end{lem}
\begin{proof}
    Available in Appendix~\ref{sect:app:lem_exist}.
\end{proof}

Lemma~\ref{lem:exist} describes a power flow solution existence condition which states that operating points exist for all $p_\ell \in [\ubar{p}_\ell,\bar{p}_\ell]$ as long as an operating point exists for the high-loading condition. Further, this lemma shows that the high-voltage solutions for any $p_\ell \in \mathcal{P}$ are bounded by those for $\ubar{p}_\ell$ and $\bar{p}_\ell$. The next result derives a simple certificate to check if the power flow solution is the high-voltage one under nonnegative power injections.

To facilitate subsequent discussion, we rearrange~\eqref{eq:pfe_ori} into the following fixed-point form:
\begin{equation}
    V_\ell = G_{p_\ell}(V_\ell) \triangleq E + Z_{\ell \ell}\diag\{p_{\ell}\}r(V_\ell),\label{eq:pfe_org}
\end{equation}
where we denote $Z_{\ell \ell} = Y^{-1}_{\ell \ell}$, $E = -Y_{\ell\ell}^{-1}Y_{\ell s} V^{\mathrm{ref}}$, and $r(V_\ell) = [1/V_{\ell 1},\cdots,1/V_{\ell n_\ell}]^\top$ yields the element-wise reciprocal of vector $V_\ell$.

\begin{prop}\label{prop:bound}
    Let the power injection $p_\ell = \bar{p}_\ell \ge \mathbbold{0}$ and source voltage $V^\mathrm{ref} > \mathbbold{0}$ be given. When
    \begin{equation}
        ||Z_{\ell\ell} \diag^{-2}\{E\} \bar{p}_{\ell}||_\infty < 1, \nonumber
    \end{equation}
    there exists a unique solution to \eqref{eq:pfe_ori} in $[E, +\infty)$.
\end{prop}
\begin{proof}
    Available in Appendix \ref{sect:app:prop_bound}.
\end{proof}

Lemma~\ref{lem:exist} provides the following two implications that align with engineering observations:
\begin{itemize}
    \item[1.] (\emph{``Monotonicity'' with respect to solution existence}) We only need to verify that the system is solvable at the high-loading condition $\ubar{p}_\ell$ to ascertain the existence of power flow solution at any loading condition $p_\ell \ge \ubar{p}_\ell$. 
    \item[2.] (\emph{Monotonicity with respect to voltage}) The high- and low-loading solutions jointly define solution bounds for power flow solutions under all loading conditions in between, where the high-loading solution provides the lower bound and low-loading solution provides the upper bound.
\end{itemize}
These implications help reduce computational efforts in DC network operations. For example, one only needs to examine whether the system has an operating point in the high-loading condition to certify all the other cases. Moreover, these conditions are key for developing a solution algorithm for~\eqref{eq:prob1}.

\subsection{Robust DN-OPF}
As we have characterized a cluster of operating points with respect to given generator setpoints, the remaining task is to design the setpoints to steer the cluster into the desired set 
(i.e., $\mathcal{V}^s_\ell \cap \mathcal{V}_\ell$)
and to reduce system operational costs.  

To accomplish these two objectives, an optimization problem can be formulated as
\begin{subequations}\label{eq:prob_final_1}
	\begin{align}
	\noalign{\textbf{R. } \textbf{DN-OPF*:}}
	\min_{\substack{V^{\mathrm{ref}}\in \mathcal{V}^{\mathrm{ref}},  s>0}} \quad & f(V^{\mathrm{ref}},V_\ell^*), \\
	\text{s.t.} \quad & \ubar{p}_\ell=\diag\{V_{\ell-}\}\left( Y_{\ell \ell}\,V_{\ell-}+Y_{\ell s}\,V^{\mathrm{ref}}\right),\label{eq:con_pfe_low} \\
	& \bar{p}_\ell=\diag\{V_{\ell+}\}\left( Y_{\ell \ell}\,V_{\ell+}+Y_{\ell s}\,V^{\mathrm{ref}}\right),\label{eq:con_pfe_up}\\
	& E = -Z_{\ell\ell}\, Y_{\ell s}\, V^{\mathrm{ref}}, \;E \le V_{\ell+},\label{eq:con_E} \\
	& \mathbbold{1}\,s\, < \,E,\; Z_{\ell\ell}\,\bar{p}_\ell < \mathbbold{1}\,s^2, \label{eq:con_pfe_org_up}\\
	& V_{\ell-} \in \mathcal{V}^s_\ell \cap \mathcal{V}_\ell, \; V_{\ell+} \in \mathcal{V}^s_\ell \cap \mathcal{V}_\ell,\label{eq:con_range}\\
	& \eqref{eq:con_pfe_ori}, \; \eqref{eq:con_pfe_ori_2}.
	\end{align} 
\end{subequations}
Constraints~\eqref{eq:con_pfe_low} and~\eqref{eq:con_pfe_up} represent the DN power flow equations for the high- and low-loading conditions, constraints~\eqref{eq:con_E}--\eqref{eq:con_pfe_org_up} enforce the condition of Proposition~\ref{prop:bound}, and constraints~\eqref{eq:con_range} and $V^{\mathrm{ref}}\in \mathcal{V}^{\mathrm{ref}}$ ensure the operating points are steered into the desired set.

A solution of~\eqref{eq:prob_final_1} yields generator setpoints which ensure robust stability and feasibility, as stated by the following result:
\begin{theorem}\label{thm:main}
    Any solution of~\eqref{eq:prob_final_1} is a feasible point of~\eqref{eq:prob1}.
\end{theorem}
\begin{proof}
    Available in Appendix~\ref{sect:app:thm}.
\end{proof}

Theorem~\ref{thm:main} shows that the feasibility of~\eqref{eq:prob_final_1} implies the feasibility of~\eqref{eq:prob1}. Problem~\eqref{eq:prob_final_1} only contains linear and quadratic constraints whose structure resembles that of a classic \mbox{DN-OPF} problem~\eqref{eq:prob_ori}. Thus, existing DN-OPF algorithms can be leveraged to solve~\eqref{eq:prob_final_1}. We use the following algorithm to summarize the main result of the paper:
\begin{algorithm}
	\caption{Find $V^{\mathrm{ref}}$ for SIP~\eqref{eq:prob1}}\label{alg:main}
	\textbf{Input:} System matrices $A$, $B$, $C$, $D$, load uncertainty set $\mathcal{P}_\ell$, constraint sets $\mathcal{V}_\ell$ and $\mathcal{V}^{\mathrm{ref}}$.\\
	\textbf{Output:} A solution $V^{\mathrm{ref}}$.
	\begin{algorithmic}[1]
		\item[Step 1:] Solve~\eqref{eq:prob_hivo} in Appendix~\ref{sect:app:set} to find individual steady-state voltage lower bound $V^2_{\ell k-}$, for all $k \in \mathcal{N}_\ell$.
		\item[Step 2:] Construct $\Delta_0$.
		\item[Step 3:] Solve GEVP~\eqref{eq:prob_gevp} to find $\alpha$ and all $V^{\mathrm{s}}_{\ell k},k\in\mathcal{N}_\ell$.
		\item[Step 4:] Find robust stability set $\mathcal{V}^{s}_\ell$.
		\item[Step 5:] Solve problem~\eqref{eq:prob_final_1} to find $V^{\text{ref}}$.
	\end{algorithmic}
\end{algorithm}

Using Algorithm~\ref{alg:main}, we can compute generator setpoints that ensure robust stability and feasibility. As previously discussed, the algorithm can be executed efficiently with existing tools.


\begin{rem}[Actual Operating Points]
   Our work can ensure the robust stability and feasibility of the actual operating points. Like a general power system, a DC network usually operates at a high-voltage power flow solution~\cite{john2015resistive}. As discussed in Appendix~\ref{sect:app:prop_bound}, all high-voltage solutions are bounded from above by that of the low-loading condition and from below by that of the high-loading condition. Note that $V_{\ell +}$ is the high-voltage solution at the low-loading condition, hence it is a tight upper bound for the system's high-voltage solutions. In addition, $V_{\ell -}$ is a solution at the high-loading condition, which is element-wise less than or equal to the high-voltage solution. Hence, any high-voltage solution must reside in the range $[V_{\ell-},V_{\ell+}]$, which is entirely steered into the desired set. 
\end{rem}

\begin{rem}[Conservativeness]
   Numerical studies show that the proposed approach has limited conservativeness. First, as discussed in Section~\ref{sec:robstab} the stability condition has limited conservativeness. Second, numerical studies suggest that the interval bound $[V_{\ell-},V_{\ell+}]$ usually has no gap with respect to the actual operating point variation range. 
\end{rem}

\begin{rem}[Computational Tractability]
   The tasks needed for Algorithm~\ref{alg:main} are computationally tractable. As discussed in Section~\ref{sec:robstab}, Steps 2 and 4 only involve simple algebraic calculations that pose trivial computational burdens. As shown in Appendix~\ref{sect:app:set}, Step 1 can be approximated by SOCPs. Step 3 involves solving a series of SDPs. There exists computationally efficient algorithms to solve both problems as well~\cite{boyd2004convex}. Since problem~\eqref{eq:prob_final_1} resembles standard deterministic DN-OPF problems, Step~5 can be completed using existing optimization solvers like IPOPT~\cite{Wachter2006} that have been shown to be effective for solving OPF problems.
\end{rem}

\section{Case Studies}\label{sec:simu}
This section demonstrates the validity of the proposed algorithm using simulation case studies. The optimization problems are solved using IPOPT~%
\cite{Wachter2006}, and the simulations are performed in Matlab/Simulink. 
\subsection{Efficacy and Conservativeness}~\label{sec:simu_1}
We first focus on an example DC network whose topology and bus types are the same as the IEEE 14-bus system. The parameters of the DC network given in Table~\ref{table:simu1_spec} are chosen according to existing DC network case studies~%
\cite{liu2018existence,salo2007low}. 

We consider three scenarios:~1)~all load case, where power outputs of each constant power component varies in~$[-50~\text{kW},0]$;~2)~all generation case, where power outputs of each constant power component varies in~$[0,50~\text{kW}]$; and~3)~all mixed case, where each constant power component varies in~$[-50~\text{kW},50~\text{kW}]$.

We use the first scenario to demonstrate that the generator setpoints designed using our results ensure robust stability and feasibility. For this case study, we impose operational bounds of $[450~\text{V}, 550~\text{V}]$ on the generator and CPL voltages, which allows a $\pm 0.1$ p.u. deviation when $500~\text{V}$ is set as base voltage. The objective function minimizes the generation costs at the nominal operating condition, which is set to be $25$~kW.

\begin{table}[t]
	\caption{Parameters for the 14-bus DC network case study} 
	\label{table:simu1_spec}
	\begin{center}
		\begin{tabular}{cccccc}
			\hline\hline
			$R_{sk}$ & 0.05 $\Omega$ & $R_{lj}$ & 5 $\Omega$ & $R_{\mathrm{c}p}$ & 0.05 $\Omega$ \\ 
			$L_{\mathrm{c}p}$ & 3mH & $C_{sk}$ & 0.75mF & $C_{lj}$ & 0.9mF \\
			\hline\hline
		\end{tabular}
	\end{center}
	
\end{table}

\begin{figure}[!t]
	\centering
	\includegraphics[width=0.4\textwidth]{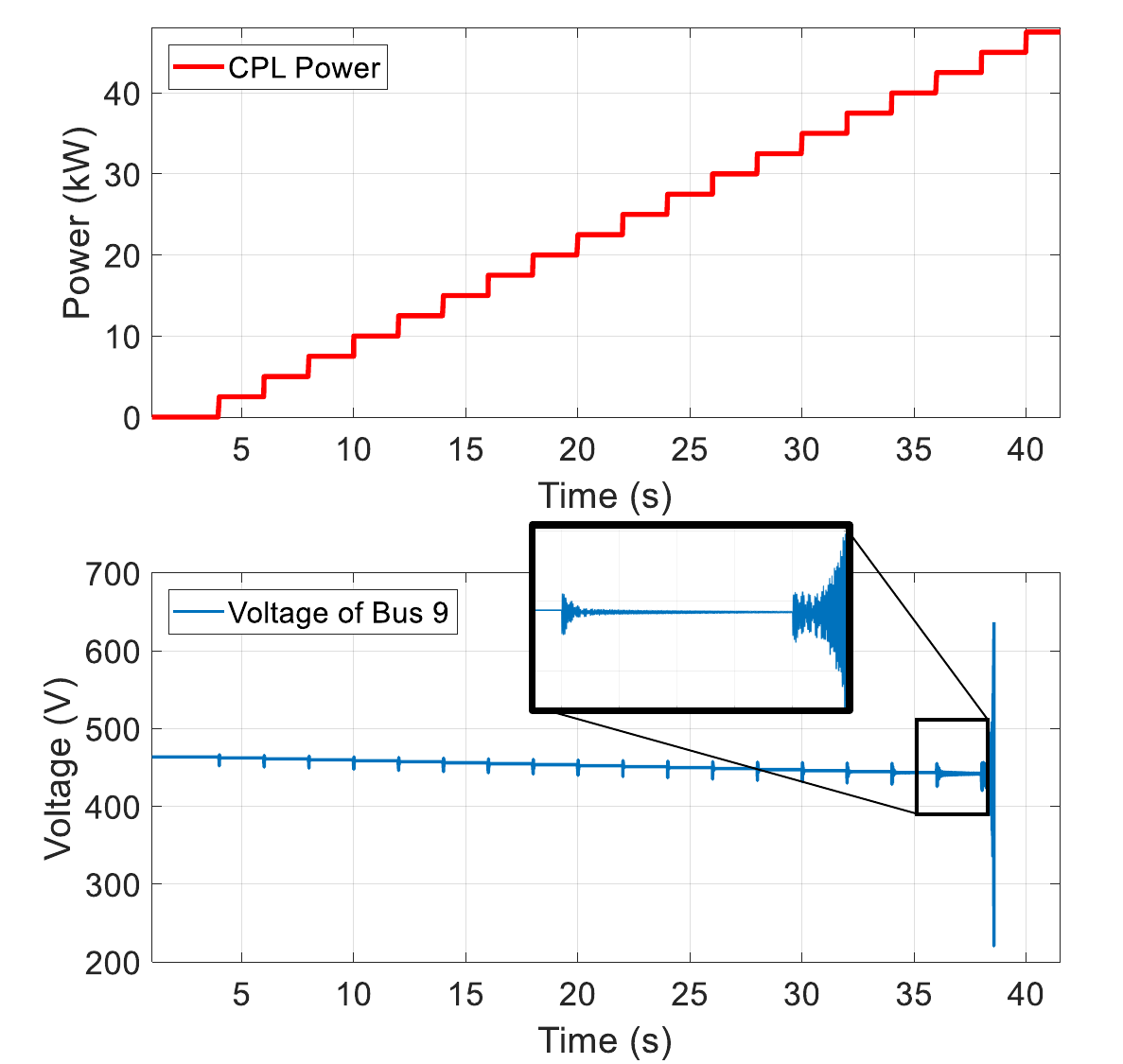}
	\caption{\label{fig:case1_uns}DC network instability if only the nominal loading condition is considered.}
\end{figure}

When we ignore the range of possible uncertainty realizations, the solution to the DN-OPF problem~\eqref{eq:prob_ori} yields setpoints of the five generators as $[481.8, 489.7, 481.2, 480.6, 486.5]$~V. We apply these setpoints and consider a uniform increase in load demands of $2.5$~kW every $2.5$ seconds. As shown in Fig.~\ref{fig:case1_uns}, the system becomes unstable at approximately $40$ seconds when the loads are around $45$~kW. In the zoomed-in figure, one can see that by approximately 38 seconds, the system is already in an oscillatory state. This shows the need to consider stability properties for operating point design, especially in systems with significant uncertainties.

\begin{figure}[!t]
	\centering
	\includegraphics[width=0.4\textwidth]{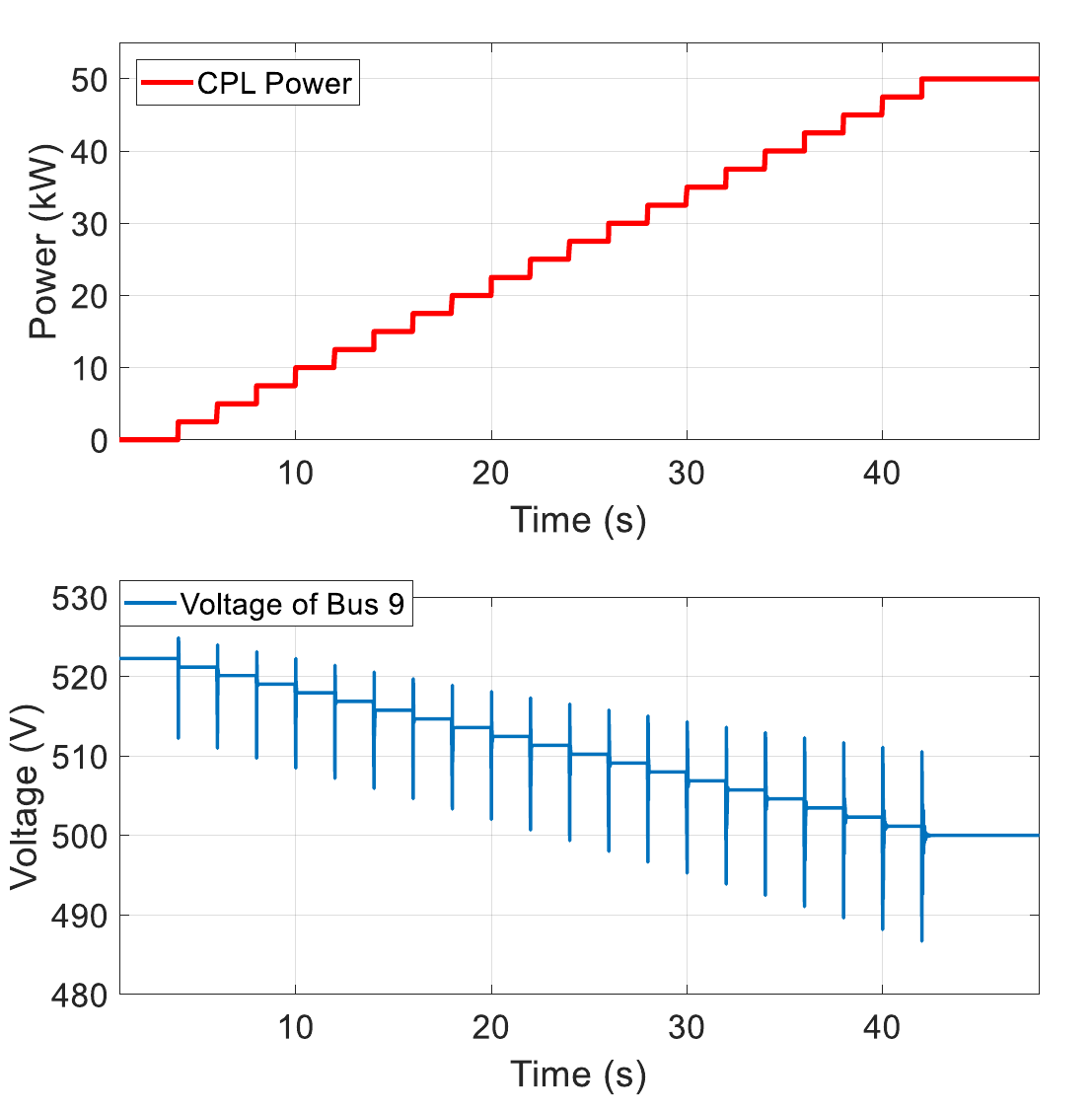}
	\caption{\label{fig:case1_stb}DC network operation is stable for all conditions using the proposed algorithm.}
\end{figure}

\begin{figure}[!t]
	\centering
	\includegraphics[width=0.4\textwidth]{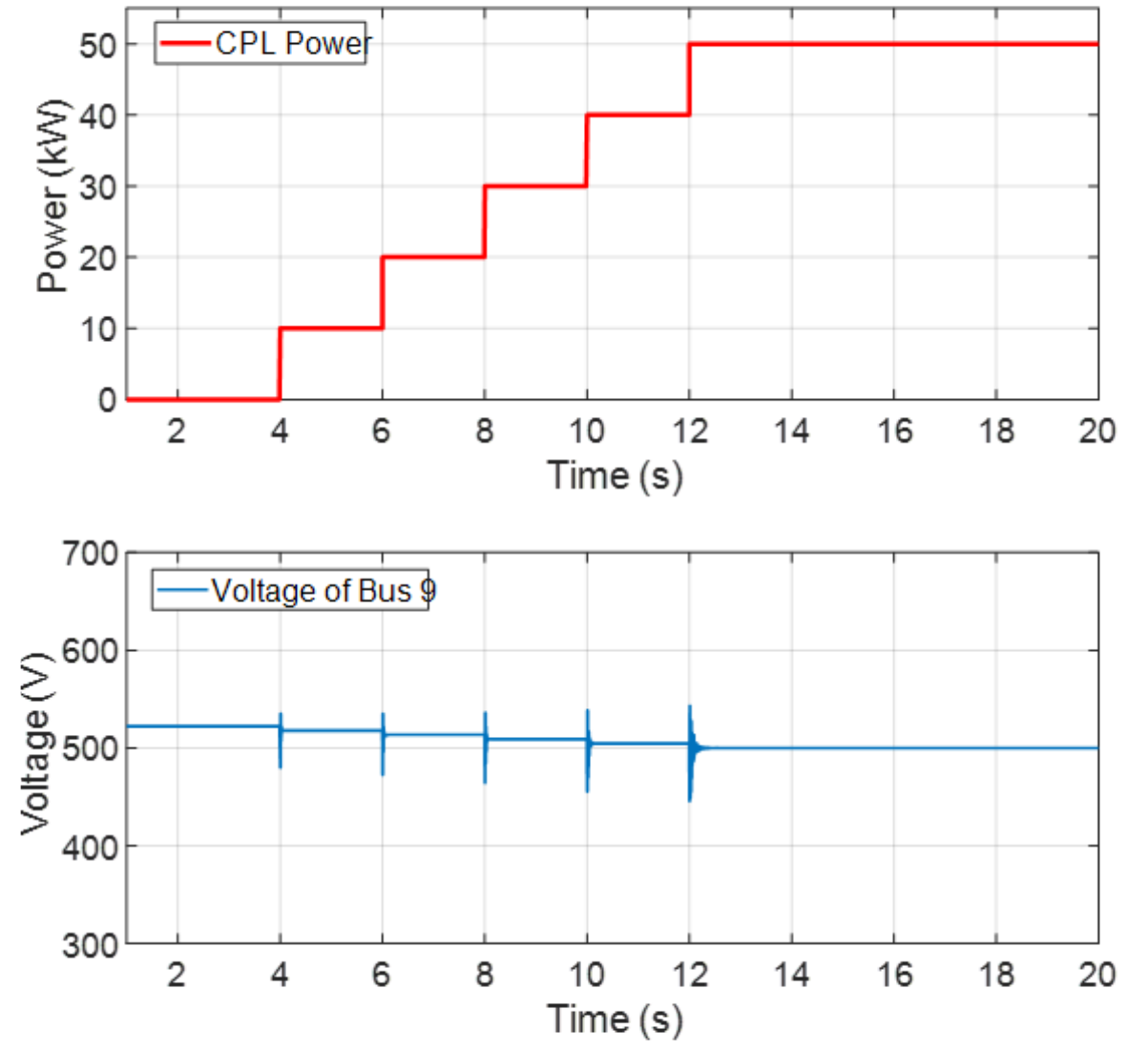}
	\caption{\label{fig:rev1_step}DC microgrid remains stable with large load variations}
\end{figure}

In comparison, we formulate the optimization problem~\eqref{eq:prob_final_1} using our algorithm. Applying the stability analysis approach developed in Section~\ref{sec:robstab} shows that the system is always robustly stable if the steady-state CPL voltage is higher than $500$~V. Problem~\eqref{eq:prob_final_1} yields the following setpoints: $[543.5, 550.0, 542.8, 542.1, 549.3]$~V. Using these setpoints results in robust stability for the entire range of load demands, as empirically corroborated by Fig.~\ref{fig:case1_stb}. Observe that the system remains stable when the loads are increased at the same rate as in the previous test. In addition, we test the stability of the system after large load step changes, where all the loads increase from $0$~kW to $50$~kW in five steps. The simulation results are shown in Fig.~\ref{fig:rev1_step}. The system is still stable despite the large load increases.

\color{black}
Moreover, we show that our approach can recover the exact power flow solution variation range. Table~\ref{table:est. error} presents the difference of the obtained $[V_{\ell -},V_{\ell +}]$ relative to the actual lower and upper bounds observed in the simulation for the three scenarios. 
Observe that we find the exact variation range of the system's actual operating points. This demonstrates the limited conservativeness of the proposed algorithm.
\begin{table}[!t]
	\caption{Recovering Exact Operating Point Variation Range}
	\label{table:est. error}
	\begin{center}
		\begin{tabular}{l|llll}
			\hline
			& All gen.  & All loads & All mixed\\ 
			\hline
			Gap for upper bound & 0  & 0  & 0  \\ 
			
			Gap for lower bound & 0 & 0 & 0  \\
			\hline
		\end{tabular}
	\end{center}
\end{table}

\begin{table}[!t]
	\caption{Comparison of computation time (in sec.) for solving~\eqref{eq:prob_ori} and~\eqref{eq:prob_final_1}}
	\label{table:compu_compa}
	\begin{center}
		\begin{tabular}{l|llllllll}
			\hline
			& 9-bus  & 39-bus & 118-bus & 300-bus & 2383-bus\\ 
			\hline
			DN-OPF~\eqref{eq:prob_ori} & $0.04$ & $0.09$ & $0.17$ & $0.26$ & $0.87$ \\ 
			Problem~\eqref{eq:prob_final_1} & $0.08$ & $0.16$ & $0.50$ & $1.77$ & $196.21$ \\
			\hline
		\end{tabular}
	\end{center}
\end{table}
\subsection*{B. Computational Efficiency}
We have tested the computational tractability of DN-OPF problem~\eqref{eq:prob_final_1} on DC networks with the same topology and bus types as the IEEE~\mbox{9-,}~\mbox{39-},~\mbox{118-},~\mbox{300-}, and~\mbox{2383-bus} systems. To summarize the results, Table~\ref{table:compu_compa} compares the average CPU time in IPOPT for solving problem~\eqref{eq:prob_final_1} and the traditional DN-OPF problem~\eqref{eq:prob_ori}, averaged over 10 tests for each system. Observe that the proposed optimization problem has a similar computational complexity as the traditional \mbox{DN-OPF} problem for systems with moderate sizes, and is still reasonably tractable for large-scale systems like 2383-bus system. This verifies the tractability of our algorithm.
\section{Conclusion}\label{sec:sum}
\color{black}This paper has developed an algorithm for solving stability-constrained OPF problems in DC networks under uncertainty. Such problems are usually intractable due to infinitely many constraints. Our algorithm uses computationally efficient approaches to transform the problem into a tractable counterpart that resembles a traditional DN-OPF problem such that existing tools can be employed. We first derive a robust stability set within which any operating point is guaranteed to be robustly stable. We then use a sufficient condition which ensures the existence of feasible operating points in this set for all uncertainty realizations in the specified uncertainty set. Low conservativeness and high computational efficiency of the proposed algorithm are demonstrated using various test cases. In future research, we will investigate the application of the algorithm to DN-OPF problems with contingency constraints.

\appendices
\section{Initial Guess for Robust Stability Set}~\label{sect:app:set}
The suggested set is motivated by the outer convex approximation of an OPF problem. The main idea is to find an outer approximation of the domain of $\delta$. To find a tight approximation, we look into the coupling between $p_{\ell}$ and $V^2_{\ell}$ since they constitute $\delta$. For different power flow solutions and power profiles, the value of $\delta$ varies. Let $\delta_{k+}$ and $\delta_{k-}$ be the upper and lower boundaries for $\delta_k$. In addition, we use $V^2_{\ell k+}$ and $V^2_{\ell k-}$ to represent a pair of upper and lower bounds of $V^2_\ell$ for any feasible solution. Recall that each power injection ranges from nonnegative to non-positive, hence $\delta_{\ell j+} = -\ubar{p}_{\ell j} / V_{\ell j-}^2$ and $\delta_{\ell j-} = -\bar{p}_{\ell j} / V_{\ell j-}^2$.

We only need to find $V^2_{\ell k-}$ in order to obtain $\delta_{k+}$ and $\delta_{k-}$. Such a bound can be found through solving the following problem:
\begin{equation}\label{eq:prob_hivo}
    \hspace{-0.02in}\textbf{Lower Bound:}\, \min_{p_\ell \in \mathcal{P}_\ell, V^{\mathrm{ref}}\in \mathcal{V}^{\mathrm{ref}}} V_{\ell k}^2,\, \text{s. t. }~\eqref{eq:pfe_ori},\, V_\ell \in \mathcal{V}_\ell.
\end{equation}
Notice that~\eqref{eq:prob_hivo} is a conventional optimal power flow problem with a QCQP formulation. Even though there exist efficient numerical solvers for such problems, the nonconvexity of~\eqref{eq:prob_hivo} may still raise concerns regarding the solution of $n_\ell$ of these problems. Since we only need to find a lower bound on the steady-state voltage, one could instead apply a second-order cone programming (SOCP) relaxation method to find a lower bound of the global optimum of~\eqref{eq:prob_hivo}.\footnote{With certain conditions on the problem structure, the SOCP relaxation of~\eqref{eq:prob_hivo} is exact. Detailed discussion can be found in~\cite{gan2014optimal,li2018optimal}.} 

\section{Proof of Proposition~\ref{prop:stab}}\label{sect:app:prop_stab}
\begin{proof}
    From the definition of $\Delta$, we have $\ubar{\delta}_k < 0 < \bar{\delta}_k$ for all $k = 1,\ldots,n_\ell$. Suppose that $\alpha > 0$ and $P \succ 0$ are solutions of GEVP~\eqref{eq:prob_gevp}. From Lemma~\ref{lem:LMIs}, any $\delta \in \alpha \Delta$ satisfies $PJ(\delta)+J^\top(\delta)P \preceq 0$. 
    
    Suppose $V_\ell$ and $p_\ell$ are arbitrary elements of $\mathcal{V}^\mathrm{s}_\ell$ and $\mathcal{P}_\ell$, respectively.
    With a slight abuse of notation, let $\delta = [-p_{\ell 1}/V^2_{\ell 1},\cdots, -p_{\ell n_\ell}/V^2_{\ell n_\ell}]^\top$. It suffices to demonstrate that $\delta \in \alpha \Delta$ to show the validity of the proposition. For each entry of $\delta$, there are three types of cases to consider: 1)~$p_{\ell k} = 0$; 2)~$p_{\ell k} > 0$; and 3)~$p_{\ell k} < 0$. For case~1), the proof is trivial as $\delta_k = 0$ and it must lie in $[\alpha \ubar{\delta}_k,\alpha \bar{\delta}_k]$. For case~2), we have the following derivation, $\alpha \ubar{\delta}_k \leq  \frac{-\bar{p}_{\ell k}}{(V^\mathrm{s}_{\ell k})^2} \leq \delta_k = \frac{-p_{\ell k}}{V^2_{\ell k}} < 0 < \alpha \bar{\delta}_k$,
    where we have used 1) $(V^\mathrm{s}_{\ell k})^2 \geq -\bar{p}_{\ell k}/(\alpha \ubar{\delta}_k)$; 2) $\ubar{\delta}_k < 0 < \bar{\delta}_k$; 3) $0 < p_{\ell k} \leq \bar{p}_{\ell k}$; and 4) $0 < (V^\mathrm{s}_{\ell k})^2 \leq V^2_{\ell k}$. The proof in regards to case 3) is similar and is omitted for brevity.
\end{proof}
\section{Proof of Lemma~\ref{lem:exist}}
\label{sect:app:lem_exist}
\begin{proof}
    To facilitate our discussion, we write $Y_{\ell\ell} = Y_{\ell\ell}^\mathrm{d} + Y_{\ell\ell}^\mathrm{off}$ where $Y_{\ell\ell}^\mathrm{d}$ and  $Y_{\ell\ell}^\mathrm{off}$ are the diagonal and off-diagonal parts of $Y_{\ell\ell}$. The power flow equation \eqref{eq:pfe_ori} can then be rearranged as
    \begin{equation} \label{eq:pfe_int}
        \diag\{V_\ell\}Y_{\ell\ell}^\mathrm{d}V_\ell = p_\ell - \diag\{V_\ell\}Y_{\ell\ell}^\mathrm{off}V_\ell - \diag\{V_\ell\}Y_{\ell s} V^\mathrm{ref}.
    \end{equation}
    Since $Y_{\ell\ell}^\mathrm{d} \succ 0$, we can left-multiply both sides of \eqref{eq:pfe_int} by $\left(Y_{\ell\ell}^\mathrm{d}\right)^{-1}$ to obtain the following reformulation of \eqref{eq:pfe_ori}:
    \begin{multline} \label{eq:pfe_int2}
        \diag\{V_\ell\}V_\ell = \left(Y_{\ell\ell}^\mathrm{d}\right)^{-1} \big( p_\ell - \diag\{V_\ell\}Y_{\ell\ell}^\mathrm{off}V_\ell \\ - \diag\{V_\ell\}Y_{\ell s} V^\mathrm{ref} \big).
    \end{multline}
    Define $W \triangleq \diag\{V_\ell\}V_\ell$, then \eqref{eq:pfe_int2} represents a fixed-point mapping of $W$ which can be written more compactly as
    \begin{equation} \label{eq:fppf}
     W = F_{p_{\ell}}(W). 
    \end{equation}
    
    We know from \cite[Lemma 2]{dvijotham2017high} that for any $p_\ell$ there exists a $W_{\max}$ such that $W_{\max} \ge W$ for any $W$ satisfying $W \le F_{p_{\ell}}(W)$. In addition, we have $W_{\max} \ge F_{p_{\ell}}(W_{\max})$. Since $(Y_{\ell\ell}^{\mathrm{d}})^{-1}, -Y_{\ell\ell}^{\mathrm{off}}, - Y_{\ell s}$ are all positive matrices, $F_{p_\ell}(\cdot)$ is an increasing function: $X \ge Y$ implies $F_{p_\ell}(X) \ge F_{p_\ell}(Y)$. Similarly, $F_{p_\ell}(W)$ is an increasing function with respect to $p_\ell$. Let $V_{\ell -}$ be the high-voltage solution for $p_\ell = \ubar{p}_\ell$ and denote $\diag\{V_{\ell -}\}V_{\ell -}$ by $W_-$. Let $p_\ell \ge \ubar{p}_\ell$ be given. By monotonicity of $F$ with respect to $p_\ell$ we know $F_{p_\ell}(W_-) \ge F_{\ubar{p}_\ell}(W_-) = W_-$. By definition of $W_{\max}$, we have $W_{\max} \ge F_{p_{\ell}}(W_{\max})$. Furthermore, $F_{p_\ell}(W_-) \le F_{p_\ell}(W) \le F_{p_\ell}(W_{\max})$ holds for any $W\in[W_-, W_{\max}]$ by monotonicity of $F_{p_\ell}(W)$. It follows from the inequalities above that $W_- \le F_{p_\ell}(W) \le W_{\max}$ for any $W\in[W_-, W_{\max}]$. In other words, $[W_-, W_{\max}]$ is an invariant set for \eqref{eq:fppf}. It follows from Brouwer's fixed-point theorem that \eqref{eq:fppf} admits a fixed point in $[W_-, W_{\max}]$. By definition of $W_{\max}$, the high-voltage solution lies in the set too. We have thus shown that the high-voltage solution exists for all $p_\ell \in \mathcal{P}$ as long as $V_{\ell -}$ exists and it is increasing with respect to load power $p_\ell$. The last statement follows. 
\end{proof}

\section{Proof of Proposition~\ref{prop:bound}}\label{sect:app:prop_bound}
\begin{proof}
     Banach's fixed-point theorem is used to establish the proof. It states that for any contraction mapping $G$ on a complete metric space mapping a set to itself, there is a unique point $x_0$ in the set such that $G(x_0) = x_0$.
     
     We show that $G_{\bar{p}_\ell}(\cdot)$ is a self map on $[E,+\infty)$. Notice that $E>\mathbbold{0}$, and $G_{\bar{p}_\ell}(\cdot) = E + Z_{\ell \ell}\diag\{\bar{p}_\ell\} r(V_\ell)$ with positive matrix $Z_{\ell \ell}$ and nonnegative $\bar{p}_\ell$. Hence, the second term on the right-hand side must be nonnegative. Thus, for all $V_\ell \in [E,+\infty)$ we have $G_{\bar{p}_\ell}(V_\ell) \ge E$. 
     
     Next, we show that $G_{\bar{p}_\ell}(\cdot)$ is a contraction mapping on $[E,+\infty)$. For any $V^a, V^b \in [E,+\infty)$ , the infinity norm of $G_{\bar{p}_\ell}(V^a) - G_{\bar{p}_\ell}(V^b)$ is upper bounded as follows:
     \begin{align}
         & \hspace{0.18in} ||G_{\bar{p}_\ell}(V^a) - G_{\bar{p}_\ell}(V^b)||_\infty \nonumber\\
         &= ||Z_{\ell \ell} \diag\{\bar{p}_\ell\} \left( r(V^a) - r(V^b) \right) ||_\infty \nonumber\\
         &\leq ||Z_{\ell\ell} \diag^{-2}\{E\} \bar{p}_{\ell}||_\infty \cdot ||V^a - V^b||_\infty, \label{eq:prop_bound_E}
     \end{align}
     where $||Z_{\ell\ell} \diag^{-2}\{E\} \bar{p}_{\ell}||_\infty$ is a positive constant less than 1. Hence, $G_{\bar{p}_\ell}(\cdot)$ is a contraction mapping on $[E, +\infty)$. 
     
     From the Banach fixed-point theorem, $G_{\bar{p}_\ell}(\cdot)$ has a unique fixed point in $[E,+\infty)$, which is clearly the high-voltage solution.
\end{proof}

\balance

\section{Proof of Theorem~\ref{thm:main}} \label{sect:app:thm}
\begin{proof}
     Suppose $V^{\mathrm{ref}}$, $V_{\ell -}$, $V_{\ell +}$, and $s$ are solutions of~\eqref{eq:prob_final_1}. We only need to show $V^{\mathrm{ref}}$ ensures robust stability and robust feasibility as defined in Definition~\ref{def:feas+stab}.
    
    First, when~\eqref{eq:prob_final_1} is feasible, the condition of Proposition~\ref{prop:bound} is satisfied with infinity-norm since constraints~\eqref{eq:con_pfe_low}--\eqref{eq:con_pfe_org_up} hold. Lemma~\ref{lem:exist} says that there must exist a power flow solution $V_{\ell}$ in $[V_{\ell -},\, V_{\ell +}]$ for any $p_\ell \in [\ubar{p}_\ell,\, \bar{p}_\ell]$.
    
    Second, from~\eqref{eq:con_range} the interval $[V_{\ell -},\, V_{\ell +}]$ belongs to $\mathcal{V}^s_\ell \cap \mathcal{V}_\ell$. Hence, the power flow solution $V_{\ell} \in \mathcal{V}^s_\ell \cap \mathcal{V}_\ell$ as well. Thus, the system power flow solution satisfies both the operational constraints and the stability constraints. 
    
    Since $V_{\ell} \in \mathcal{V}_{\ell}$ and $V^{\mathrm{ref}} \in \mathcal{V}^{\mathrm{ref}}$, the system is robustly feasible. Since $V_{\ell} \in \mathcal{V}^s_{\ell}$, by Proposition~\ref{prop:stab} the system is also robustly stable. Thus, the proof is complete.
\end{proof}

\bibliographystyle{IEEEtran}
\bibliography{DCOPF}

%
%
%
%
%

\end{document}

%% file: nomenclature.tex
\printnomenclature
\nomenclature[P]{$n$}{Total number of components}
\nomenclature[P]{$n_s$}{Total number of generators}
\nomenclature[P]{$n_\ell$}{Total number of loads}
\nomenclature[P]{$n_{\mathrm{c}}$}{Total number of power lines}
\nomenclature[P]{$p_{\ell}$}{Power injection of loads}
\nomenclature[P]{$p^*_{\ell}$}{Nominal power injection of loads}
\nomenclature[P]{$R_{\ell}$}{Resistance for loads}
\nomenclature[P]{$C_{\ell}$}{Capacitance for loads}
\nomenclature[P]{$R_{s}$}{Resistance for sources}
\nomenclature[P]{$C_{s}$}{Capacitance for sources}
\nomenclature[P]{$R_{\mathrm{c}}$}{Resistance for power lines}
\nomenclature[P]{$L_{\mathrm{c}}$}{Inductance for power lines}
\nomenclature[P]{$Y_{\ell \ell}$}{Admittance submatrix for load connectivity}
\nomenclature[P]{$Y_{\ell s}$}{Admittance submatrix for load-generator connectivity}
\nomenclature[p]{$\delta$}{Composite uncertainty}
\nomenclature[V]{$i_{o}(t)$}{Current flowing into power line}
\nomenclature[V]{$i_{d}(t)$}{Current flowing out of power line}
\nomenclature[V]{$V^{\mathrm{ref}}$}{Voltage set-point for generators}
\nomenclature[V]{$v_{s}(t)$}{External voltage of generators}
\nomenclature[V]{$v_{\ell}(t)$}{Voltage across capacitor of loads}
\nomenclature[V]{$i_{\mathrm{c}}(t)$}{Current through power lines}
\nomenclature[V]{$x$}{System states}
\nomenclature[V]{$x^e$}{System equilibrium}
\nomenclature[V]{$I_{\mathrm{c}}$}{Power line current at equilibrium}
\nomenclature[V]{$V_s$}{Generator terminal voltage at equilibrium}
\nomenclature[V]{$V_\ell$}{Load terminal voltage at equilibrium}
\nomenclature[S]{$\mathcal{N}_s$}{Index set of generators}
\nomenclature[S]{$\mathcal{N}_\ell$}{Index set of loads}
\nomenclature[S]{$\mathcal{E}_{\mathrm{c}}$}{Index set of power lines}
\nomenclature[S]{$\mathcal{P}_\ell$}{Set of load profiles}
\nomenclature[S]{$\mathcal{V}^{\mathrm{ref}}$}{Generator voltage set-point constraint set}
\nomenclature[S]{$\mathcal{V}_\ell$}{Generator voltage set-point constraint set}
\nomenclature[S]{$\Delta$}{Perturbation set}